\documentclass[10pt]{amsart}
     \makeatletter
     \def\section{\@startsection{section}{1}%
     \z@{.7\linespacing\@plus\linespacing}{.5\linespacing}%
     {\bfseries
     \centering
     }}
     \def\@secnumfont{\bfseries}
     \makeatother
\newtheorem{theorem}{Theorem}[section]
\newtheorem{lemma}[theorem]{Lemma}

\theoremstyle{definition}
\newtheorem{definition}[theorem]{Definition}

\theoremstyle{remark}

\numberwithin{equation}{section}
\usepackage{amsmath,amsthm,amssymb,amsbsy}

\usepackage[all]{xy}
\usepackage{chngcntr}

\counterwithin{figure}{section}

\newcommand{\MM}{\mathbb{M}}

\newcommand{\RR}{\mathbb{R}}

\newcommand{\FFF}{\mathcal{F}}

\newcommand{\df}[1]{\,\mathrm{d}#1}                         

\begin{document}

\setlength{\parindent}{0cm}
\setlength{\parskip}{0.5cm}

\title{Intervention in Ornstein-Uhlenbeck SDEs}

\author{Alexander Sokol}

\address{Alexander Sokol: Institute of Mathematics, University of
  Copenhagen, 2100 Copenhagen, Denmark}
\email{alexander@math.ku.dk}
\urladdr{http://www.math.ku.dk/$\sim$alexander}



\begin{abstract}
We introduce a notion of intervention for stochastic differential
equations and a corresponding causal interpretation. For the case of
the Ornstein-Uhlenbeck SDE, we show that the SDE resulting from a
simple type of intervention again is an Ornstein-Uhlenbeck SDE. We
discuss criteria for the existence of a stationary distribution for
the solution to the intervened SDE. We illustrate the effect of
interventions by calculating the mean and variance in the stationary
distribution of an intervened process in a particularly simple case.
\end{abstract}

\maketitle

\noindent

\section{Introduction}


Causal inference for continuous-time processes is a field in ongoing
development. Similar to causal inference for graphical models, see
\cite{MR2548166}, one of the primary objectives for causal inference for continuous-time
processes is to identify the effect of an intervention given
assumptions on the distribution and causal structure of the observed
continuous-time process.

Several flavours of causal inference are available for continuous-time
processes, see for example \cite{MR1403234,MR2575938,MR2811860}. In this
paper, we outline a causal interpretation of stochastic differential
equations and a corresponding notion of intervention, we calculate the distribution of an intervened
Ornstein-Uhlenbeck SDE, and we calculate analytical expressions
for the mean and variance of the stationary distribution of the
resulting process for particular examples of interventions.

\section{Causal interpretation of stochastic differential equations}

\label{sec:CausalInterpretation}

Consider a filtered probability space
$(\Omega,\FFF,(\FFF_t)_{t\ge0},P)$ satisfying the usual conditions,
see \cite{MR2273672} for the definition of this and other notions related to
continuous-time stochastic processes. Let $Z$ be a $d$-dimensional
semimartingale and assume that $a:\RR^p\to\MM(p,d)$ is a Lipschitz
mapping, where $\MM(p,d)$ denotes the space of real $p\times d$
matrices. Consider the stochastic differential equation (SDE)
\begin{align}
  X^i_t &= x^i_0 + \sum_{j=1}^d \int_0^t a_{ij}(X_{s-})\df{Z^j}_s,
  \qquad i\le p.\label{eq:MainSDE}
\end{align}
By the Lipschitz property of $a$, it holds by Theorem V.7 of \cite{MR2273672} that there exists a
pathwisely unique solution to (\ref{eq:MainSDE}). The following
definition yields a causal interpretation of (\ref{eq:MainSDE}) based
on simple substitution and inspired by ideas outlined in Section 4.1 of \cite{MR2993496}.

\begin{definition}
\label{def:SDEIntervention}
Consider some $m\le p$ and $c\in\RR$. The $(p-1)$-dimensional intervened SDE arising from the
intervention $X^m := c$ is defined to be
\begin{align}
  U^i_t &= x^i_0 + \sum_{j=1}^d \int_0^t b_{ij}(U_{s-})\df{Z^j}_s
  \textrm{ for } i\le p\textrm{ with } i\neq m, \label{eq:IntervenedSDE}
\end{align}
where $b_{ij}(y_1,\ldots,y_{m-1},y_{m+1},\ldots,y_p)=a_{ij}(y_1,\ldots,c,\ldots,y_p)$, and
the $c$ is on the $m$'th coordinate. Letting $U$ be the unique
solution to the SDE and defining $Y$ by putting $Y=(U^1,\ldots,U^{m-1},c,U^{m+1},\ldots,U^p)$, we refer to $Y$ as the intervened process and write $(X|X^m := c)$ for $Y$.
\end{definition}

By Theorem V.16 and Theorem V.5 of \cite{MR2273672}, the solutions to both (\ref{eq:MainSDE})
and (\ref{eq:IntervenedSDE}) may be approximated by the Euler schemes
for their respective SDEs. Making these approximations and applying Pearl's
notion of intervention in an appropriate sense, see \cite{MR2548166}, we may
interpret Definition \ref{def:SDEIntervention} as intervening in the
system (\ref{eq:MainSDE}) under the assumption that the driving semimartingales $Z^1,\ldots,Z^d$ are
noise processes unaffected by interventions, while the processes
$X^1,\ldots,X^p$ are affected by interventions. Note that the
operation of making an intervention takes a $p$-dimensional SDE as its input and
yields a $(p-1)$-dimensional SDE as its output, and this operation is crucially dependent
on the coefficients in the SDE: These coefficients in a sense
corresponds to the directed acyclic graphs of \cite{MR2548166}. A major benefit of
causality in systems such as (\ref{eq:MainSDE}) as compared to the theory
of \cite{MR2548166} is the ability to capture feedback systems and interventions in
such feedback systems.

As the solutions to (\ref{eq:MainSDE}) and (\ref{eq:IntervenedSDE}) are defined on the same
probability space, we may even consider the process $Y - X$, where $Y=(X|X^m :=
c)$, allowing us to calculate for example the variance of the effect
of the intervention. As $Y$ and $X$ are never observed simultaneously
in practice, however, we will concentrate on analyzing the differences
between the laws of $Y$ and $X$ separately.

\section{Intervention in Ornstein-Uhlenbeck SDEs}

\label{sec:OUIntervention}

Recall that for an $\FFF_0$ measurable variable $X_0$ and for $A\in\RR^p$, $B\in\MM(p,p)$ and $\sigma\in\MM(p,d)$,
the Ornstein-Uhlenbeck SDE with initial value
$X_0$, mean reversion level $A$, mean reversion speed $B$, diffusion
matrix $\sigma$ and $d$-dimensional driving noise is
\begin{align}
  X_t &= X_0 + \int_0^t B(X_s-A)\df{s}+\sigma W_t,\label{eq:OUSDE}
\end{align}
where $W$ is a $d$-dimensional $(\FFF_t)$ Brownian motion, see Section
II.72 of \cite{MR1796539}. The unique solution to this equation is
\begin{align}
 X_t =&\exp(tB)\left(X_0-\int_0^t \exp(-sB)BA\df{s}+\int_0^t\exp(-sB)\sigma\df{W}_s\right)
\end{align}
where the matrix exponential is defined by $\exp(A)=\sum_{n=0}^\infty
A^n/n!$. This is a Gaussian homogeneous
Markov process with continuous sample paths. The following lemma shows
that making an intervention in an Ornstein-Uhlenbeck SDE yields an
SDE whose nontrivial coordinates solve another Ornstein-Uhlenbeck SDE.


\begin{lemma}
\label{lemma:SimpleOUIntervention}
Consider the Ornstein-Uhlenbeck SDE (\ref{eq:OUSDE}) with initial
value $x_0$. Fix $m\le p$ and $c\in\RR$, and let $X$ be the unique
solution to (\ref{eq:OUSDE}). Furthermore, let $Y=(X|X^m:=c)$ and let $Y^{-m}$ be the $p-1$ dimensional
process obtained by removing the $m$'th coordinate from $Y$. Let
$\tilde{B}$ be the submatrix of $B$ obtained by removing the $m$'th
row and column of $B$, and assume that $\tilde{B}$ is invertible. Then $Y^{-m}$ solves
\begin{align}
  Y^{-m}_t &= y_0 + \int_0^t \tilde{B}(Y^{-m}_s-\tilde{A})\df{s}+\tilde{\sigma} W_t,\label{eq:OUSDEInter}
\end{align}
where $y_0$ is obtained by removing the $m$'th coordinate
from $x_0$, $\tilde{\sigma}$ is obtained by removing the $m$'th row of
$\sigma$ and $\tilde{A}=\alpha-\tilde{B}^{-1}\beta$, where $\alpha$
and $\beta$ are obtained by removing the $m$'th coordinate from $A$
and from the vector whose $i$'th component is $b_{im}(c-a_m)$,
respectively, where $b_{im}$ is the entry corresponding to the $i$'th row
and the $m$'th column of $B$, and $a_m$ is the
$m$'th element of $A$.
\end{lemma}
\textit{Proof.}
By Definition \ref{def:SDEIntervention}, we have
\begin{align}
  Y^i_t =& y_0 + \int_0^t b_{im}(c-a_m)+\sum_{j\neq m} b_{ij}(Y^j_s-a_j)\df{s} + \sum_{j=1}^p\sigma_{ij}W^j_t
\end{align}
for $i\neq m$. Note that for any vector $y$, the
system of equations in $\tilde{a}$
\begin{align}
  b_{im}(c-a_m)+\sum_{j\neq m} b_{ij}(y_j-a_j) = \sum_{j\neq
    m}b_{ij}(y_j-\tilde{a}_j)\textrm{ for }i\neq m,
\end{align}
is equivalent to the system of equations
\begin{align}
   \sum_{j\neq m}b_{ij}\tilde{a}_j
  =\left(\sum_{j\neq m} b_{ij}a_j\right)-b_{im}(c-a_m) \textrm{ for }i\neq m.
\end{align}
Since we have assumed $\tilde{B}$ to be invertible, this system of
equations has the unique solution
$\tilde{A}=\tilde{B}^{-1}(\tilde{B}\alpha-\beta)=\alpha-\tilde{B}^{-1}\beta$. For
$i\neq m$, we therefore obtain that $Y^i_t = y_0 + \int_0^t \sum_{j\neq m}
b_{ij}(Y^j_s-\tilde{a}_j)\df{s}+ \sum_{j=1}^p \sigma_{ij}W^j_t$, proving the result.
\hfill$\Box$

Recall that a principal submatrix of a matrix is a submatrix with the
same rows and columns removed. In words, Lemma
\ref{lemma:SimpleOUIntervention} states that if a particular principal
submatrix $\tilde{B}$ of the mean reversion speed is invertible, then making the
intervention $X^m:=c$ in an Ornstein-Uhlenbeck SDE results in a
new Ornstein-Uhlenbeck SDE with mean reversion speed $\tilde{B}$ and
modified mean reversion level involving the inverse of $\tilde{B}$. Now
assume that an Ornstein-Uhlenbeck SDE is given such that the solution
has a stationary initial distribution. A natural question to ask is
what interventions will yield intervened processes where stationary
initial distributions also exist. In the following, we consider this question.

Recall that a square matrix is called stable if
its eigenvalues have negative real parts and semistable if its
eigenvalues have nonpositive real parts, see \cite{MR0480586}. Theorem 4.1
of \cite{MR0260056} yields necessary and sufficient criteria for the
existence of a stationary probability measure for the solution of
(\ref{eq:OUSDE}). One criterion is
expressed in terms of the controllability subspace of of the matrix
pair $(B,\sigma)$, which is the span of the columns in the matrices
$\sigma, B\sigma, \ldots, B^{p-1}\sigma$. In the case where $\sigma$ has full column span, meaning
that the columns of $\sigma$ span all of $\RR^p$, the controllability
subspace is all of $\RR^p$, and Theorem 4.1 of \cite{MR0260056} shows that the existence of a stationary
probability measure is equivalent to $B$ being stable. The case where
$\sigma$ is not required to have full column span is more involved. 

In the following, we will restrict our attention to Ornstein-Uhlenbeck processes with
$\sigma$ having full column span. By Theorem 4.1 of \cite{MR0260056}, it then
holds that there exists a stationary distribution if and only if $B$ is stable. Furthermore,
applying Theorem 2.4 and Theorem 2.12 of \cite{MJ}, it holds
in the affirmative case that the stationary
distribution is the normal distribution with mean $\mu$ and variance
$\Gamma$ solving $B\mu=BA$ and $\sigma\sigma^t+B\Gamma+\Gamma B^t = 0$. Note that
as $B$ is stable, zero is not an eigenvalue of $B$, thus $B$ is
invertible and $\mu=A$. Also, stability of $B$ yields that $\Gamma
=\int_0^\infty e^{sB}\sigma\sigma^te^{sB^t}\df{s}$. For the
$(p-1)$-dimensional Ornstein-Uhlenbeck process resulting from an
intervention according to Lemma \ref{lemma:SimpleOUIntervention}, the
diffusion matrix $\tilde{\sigma}$ is obtained by removing the $m$'th row of
$\sigma$. As the columns of $\sigma$ span $\RR^p$, the columns of
$\tilde{\sigma}$ span $\RR^{p-1}$. Therefore, it also holds for the
intervened process that there exists a stationary distribution if and only
if the mean reversion speed is stable. We conclude that for diffusion
matrices with full column span, the existence of stationary distributions
for both the original and the intervened SDE is determined solely by
stability of the mean reversion speed matrix $B$ and the corresponding principal submatrices.

Consider a stable matrix $B$. It then holds that if all principal
submatrices of $B$ are stable, all interventions will preseve
stability of the system. We are thus lead to the question of when a principal submatrix of a
matrix is stable. That stability does not in general lead to stability
of principal submatrices may be seen from the following example. Define $B$ by putting
\begin{displaymath}
  B = \left[\begin{array}{cc} 1 & 7 \\ -1 & -3 \end{array}\right].
\end{displaymath}
The matrix $B$ has eigenvalues $-1\pm i\sqrt{3}$ and is thus
stable, while the principal submatrix obtained by removing the second
row and second column trivially has the single eigenvalue $1$ and thus
is not stable, in fact not even semistable. Conversely, $-B$ has
eigenvalues $1\pm i\sqrt{3}$ and thus is neither stable nor
semistable, while the principal submatrix obtained by removing the second
row and second column of $-B$ is stable.

There are classes of matrices satisfying that all principal
submatrices are stable. For example, by the inclusion principle for symmetric matrices, see Theorem
4.3.15 of \cite{MR832183}, it follows that a principal submatrix of any symmetric stable matrix
again is stable. In general, though, it is difficult
to ensure that all principal submatrices are stable. However, there
are criteria ensuring that all principal submatices are
semistable. For example, Lemma 2.4 of \cite{MR1971090} shows that if $B$ is stable and sign symmetric, then all
principal submatrices of $B$ is semistable. Here, sign symmetry is a
somewhat involved matrix criterion, it does however hold that any
stable symmtric matrix also is sign symmetric. Furthermore, by Theorem 1 of
\cite{MR0480586}, either of the follow three properties are also
sufficient for having all principal submatrices being semistable:
\begin{enumerate}
\item $A-D$ is stable for all nonnegative diagonal $D$.
\item $DA$ is stable for all positive diagonal $D$.
\item There is positive diagonal $D$ such that $AD+DA^t$ is negative definite.
\end{enumerate}



\section{An example of a particular intervention}

\label{sec:Examples}

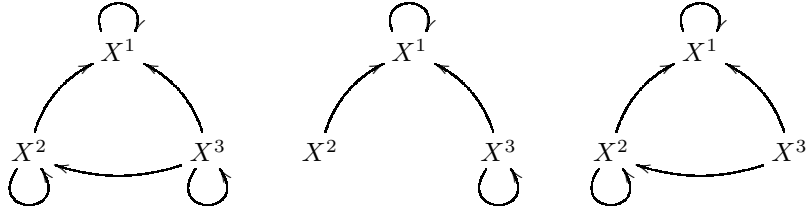
\begin{figure}[htb]
\begin{minipage}[b]{0.300\linewidth}
\vspace{0.5cm}
\centering
\begin{displaymath}
  \xymatrix@C=0.5cm{& X^1 \ar@(ul,ur)[] & \\
                    X^2 \ar@(dl,dr)[] \ar@/^8pt/[ur] & & 
                    X^3 \ar@(dl,dr)[] \ar@/_8pt/[ul] \ar@/^8pt/[ll]
                 }
\end{displaymath}
\vspace{0.5cm}
\end{minipage}
\begin{minipage}[b]{0.300\linewidth}
\vspace{0.5cm}
\centering
\begin{displaymath}
  \xymatrix@C=0.5cm{& X^1 \ar@(ul,ur)[] & \\
                    X^2 \ar@/^8pt/[ur] & & 
                    X^3 \ar@(dl,dr)[] \ar@/_8pt/[ul]
                 }
\end{displaymath}
\vspace{0.5cm}
\end{minipage}
\begin{minipage}[b]{0.300\linewidth}
\vspace{0.5cm}
\centering
\begin{displaymath}
  \xymatrix@C=0.5cm{& X^1 \ar@(ul,ur)[] & \\
                    X^2 \ar@(dl,dr)[] \ar@/^8pt/[ur] & & 
                    X^3 \ar@/_8pt/[ul] \ar@/^8pt/[ll]
                 }
\end{displaymath}
\vspace{0.5cm}
\end{minipage}
\caption{Graphical illustrations of the dependence structures of
  $(X^1,X^2,X^3)$ (left), of the dependence when making the
  intervention $X^2:=c$ (middle) and of the dependence when making the
  intervention $X^3:=c$ (right).}
\label{figure:SimpleOUDAG}
\end{figure}

Consider now a three-dimensional Ornstein-Uhlenbeck process $X$ with
$\sigma$ being the identity matrix of order three and upper diagonal mean reversion speed matrix $B$, and
assume that the diagonal elements of $B$ all are negative. As the
diagonal elements of $B$ in this case also are the eigenvalues, $B$ is then
stable, and all principal submatrices are stable as well. The interpretation of having $B$ upper diagonal is that the levels of
both $X^1$, $X^2$ and $X^3$ influence the average change in $X^1$, while only
the levels of $X^2$ and $X^3$ influence the average change in $X^2$
and only $X^3$ influences the average change in $X^3$. Figure
\ref{figure:SimpleOUDAG} illustrates this, as well as the changes to
the dependence structure obtained by making interventions
$X^2:=c$ or $X^3:=c$.

We will investigate the details of what happens to the system when
making the intervention $X^2:=c$ or $X^3:=c$. To this end, we
calculate the mean and variance in the stationary distribution for
the nontrivial coordinates in each of the intervened processes. Consider first the case of the intervention $X^2:=c$. Let $\mu$ and
$\Gamma$ denote the mean and variance in the stationary
distribution after intervention. Applying Lemma \ref{lemma:SimpleOUIntervention}, the
SDE resulting from making this intervention is a two-dimensional Ornstein-Uhlenbeck SDE
with mean reversion speed and mean reversion level
\begin{align}
  \left[\begin{array}{cc} b_{11} & b_{13} \\ 0 & b_{33} \end{array}\right]
  \quad\textrm{ and }\quad
  \left[\begin{array}{cc} a_1 \\ a_3 \end{array}\right]
  -\left[\begin{array}{cc} b_{11} & b_{13} \\ 0 & b_{33} \end{array}\right]^{-1}
  \left[\begin{array}{cc} b_{12}(c-a_2) \\ 0 \end{array}\right].
\end{align}
As we have
\begin{align}
    \left[\begin{array}{cc} b_{11} & b_{13} \\ 0 & b_{33} \end{array}\right]^{-1}
  &=\left[\begin{array}{cc} \frac{1}{b_{11}} & -\frac{b_{13}}{b_{11}b_{33}} \\ 0 & \frac{1}{b_{33}} \end{array}\right],
\end{align}
this immediately yields that
\begin{align}
    \mu &= \left[\begin{array}{cc} a_1-\frac{b_{12}}{b_{11}}(c-a_2) \\
                            a_3
     \end{array}\right].\label{eq:AsympMean}
\end{align}
As for the variance, recall that we have the representation
\begin{align}
  \Gamma &= \int_0^\infty
           \exp\left(s\left[\begin{array}{cc} b_{11} & b_{13} \\ 0 & b_{33} \end{array}\right]\right)
           \exp\left(s\left[\begin{array}{cc} b_{11} & 0 \\ b_{13} & b_{33} \end{array}\right]\right)
           \df{s}.
\end{align}
In order to calculate this integral, first consider the case
$b_{11}=b_{33}$. By Theorem 4.11 of \cite{MR2396439}, we in this case obtain
\begin{align}
  \exp\left(s\left[\begin{array}{cc} b_{11} & b_{13} \\ 0 & b_{33} \end{array}\right]\right)
  &=e^{sb_{11}}\left[\begin{array}{cc} 1 & sb_{13} \\ 0 & 1 \end{array}\right],
\end{align}
and similarly for the transpose. Applying that $\int_0^\infty x^\alpha e^{\beta x}\df{x} =
\Gamma(\alpha+1)/(-\beta)^{\alpha+1}$ for all $\alpha>-1$ and
$\beta<0$, we conclude
\begin{align}
  \Gamma
  &=  \int_0^\infty e^{2sb_{11}}
                   \left[\begin{array}{cc} 1 & sb_{13} \\ 0 & 1 \end{array}\right]
                   \left[\begin{array}{cc} 1 & 0 \\ sb_{13} & 1 \end{array}\right]\df{s}\notag\\
  &=  \int_0^\infty e^{2sb_{11}}
                   \left[\begin{array}{cc} 1+s^2b^2_{13} & 
                    sb_{13} \\ 
                    sb_{13} &
                    1 \end{array}\right]\df{s}
  =                \left[\begin{array}{cc}
                    -\frac{1}{2b_{11}}-\frac{b^2_{13}}{4b_{11}^3} & 
                    \frac{b_{13}}{4b_{11}^2} \\ 
                    \frac{b_{13}}{4b_{11}^2} &
                    -\frac{1}{2b_{11}}
                     \end{array}\right].
\end{align}
In the case $b_{11}\neq b_{33}$, we put $\zeta=b_{13}/(b_{11}-b_{33})$
and Theorem 4.11 of \cite{MR2396439} yields
\begin{align}
  \exp\left(s\left[\begin{array}{cc} b_{11} & b_{13} \\ 0 & b_{33} \end{array}\right]\right)
  &=\left[\begin{array}{cc} e^{sb_{11}} & 
                            \zeta(e^{sb_{11}}-e^{sb_{33}}) \\ 
                            0 &
                            e^{sb_{33}} \end{array}\right],
\end{align}
and we then obtain
\begin{align}
  &         \exp\left(s\left[\begin{array}{cc} b_{11} & b_{13} \\ 0 & b_{33} \end{array}\right]\right)
           \exp\left(s\left[\begin{array}{cc} b_{11} & 0 \\ b_{13} & b_{33} \end{array}\right]\right)\notag\\
  &=\left[\begin{array}{cc} e^{sb_{11}} & 
                            \zeta(e^{sb_{11}}-e^{sb_{33}}) \\ 
                            0 &
                            e^{sb_{33}} \end{array}\right]
     \left[\begin{array}{cc} e^{sb_{11}} & 
                            0 \\ 
                            \zeta(e^{sb_{11}}-e^{sb_{33}}) &
                            e^{sb_{33}} \end{array}\right]\notag\\
  &=\left[\begin{array}{cc} (1+\zeta^2)e^{2sb_{11}}
                            -2\zeta^2 e^{s(b_{11}+b_{33})}
                            +\zeta^2e^{2sb_{33}} &
                            \zeta e^{s(b_{11}+b_{33})}
                            -\zeta e^{2sb_{33}} \\
                            \zeta e^{s(b_{11}+b_{33})}
                            -\zeta e^{2sb_{33}} &
                            e^{2sb_{33}} \end{array}\right],
\end{align}
implying that
\begin{align}
  \Gamma
  &=\left[\begin{array}{cc}
     -\frac{(1+\zeta^2)}{2b_{11}}+\frac{2\zeta^2}{b_{11}+b_{33}}-\frac{\zeta^2}{2b_{33}} &
     -\frac{\zeta}{b_{11}+b_{33}}+\frac{\zeta}{2b_{33}} \\
     -\frac{\zeta}{b_{11}+b_{33}}+\frac{\zeta}{2b_{33}} &
     -\frac{1}{2b_{33}}
    \end{array}\right]\notag\\
  &=\left[\begin{array}{cc}
     -\frac{1}{2b_{11}}-\zeta^2\left(\frac{1}{2b_{11}}+\frac{2}{b_{11}+b_{33}}-\frac{1}{2b_{33}}\right) &
     \frac{\zeta(b_{11}-b_{33})}{2b_{33}(b_{11}+b_{33})} \\
     \frac{\zeta(b_{11}-b_{33})}{2b_{33}(b_{11}+b_{33})} &
     -\frac{1}{2b_{33}}
    \end{array}\right]\notag\\
  &=\left[\begin{array}{cc}
     -\frac{1}{2b_{11}}-\frac{b_{13}^2}{2b_{11}b_{33}(b_{11}+b_{33})} &
     \frac{b_{13}}{2b_{33}(b_{11}+b_{33})} \\
     \frac{b_{13}}{2b_{33}(b_{11}+b_{33})} &
     -\frac{1}{2b_{33}}
    \end{array}\right].\label{eq:AsympVar}
\end{align}
Note in particular that (\ref{eq:AsympVar}) also yields the correct
result in the case $b_{11}=b_{33}$. Next, considering the intervention
$X^3:=c$, we let $\nu$ and
$\Sigma$ denote the mean and variance in the stationary
distribution of the nontrivial coordinates after intervention. By Lemma \ref{lemma:SimpleOUIntervention}, the
result of making this interveniton is an Ornstein-Uhlenbeck SDE
with mean reversion speed and mean reversion level
\begin{align}
  \left[\begin{array}{cc} b_{11} & b_{12} \\ 0 & b_{22} \end{array}\right]
  \quad\textrm{ and }\quad
  \left[\begin{array}{cc} a_1 \\ a_2 \end{array}\right]
  -\left[\begin{array}{cc} b_{11} & b_{12} \\ 0 & b_{22} \end{array}\right]^{-1}
  \left[\begin{array}{cc} b_{13}(c-a_3) \\ b_{23}(c-a_3) \end{array}\right],
\end{align}
yielding by calculations similar to the previous case that
\begin{align}
    \nu &=\left[\begin{array}{cc}
        a_1-\left(\frac{b_{13}}{b_{11}}-\frac{b_{12}b_{23}}{b_{11}b_{22}}\right)(c-a_3) \\
                            a_2-\frac{b_{23}}{b_{22}}(c-a_3)
     \end{array}\right]\label{eq:AsympMean2}
\end{align}
and
\begin{align}
  \Sigma &=  \left[\begin{array}{cc}
     -\frac{1}{2b_{11}}-\frac{b_{12}^2}{2b_{11}b_{22}(b_{11}+b_{22})} &
     \frac{b_{12}}{2b_{22}(b_{11}+b_{22})} \\
     \frac{b_{12}}{2b_{22}(b_{11}+b_{22})} &
     -\frac{1}{2b_{22}}
    \end{array}\right].\label{eq:AsympVar2}
\end{align}
We have now calculated the mean and variance in the stationary
distribution for both intervened processes. We next take a moment to
interpret our results.

In the original system, all of $X^1$, $X^2$ and $X^3$
negatively influenced themselves, and in addition to this, $X^2$
influenced $X^1$ and $X^3$ influenced $X^1$ both directly and through
its influence on $X^2$. Based on this, we would expect that
making the intervention $X^2:=c$, the steady state of $X^3$ would not
be changed, while the steady state of $X^1$ would change, depending on
the level of influence $b_{12}$ of $X^2$ on $X^1$. This is what we see
in (\ref{eq:AsympMean}). When making the intervention $X^3:=c$,
however, we obtain a change in the steady state of $X^1$ based both on
the direct influence of $X^3$ on $X^1$, depending on $b_{13}$, but
also on the indirect influence of $X^3$ on $X^1$ through $X^2$,
depending also on $b_{23}$ and $b_{12}$. Furthermore, the steady state
of $X^2$ also changes. These results show themselves in (\ref{eq:AsympMean2}).
 
As for the steady state variance, the changes resulting from interventions are in
both cases of the same type, yielding moderately complicated
analytical expressions, both independent of $c$. This implies that
while we in most cases will be able to obtain any steady state mean
for, say, $X^1$, by picking $c$ suitably, the steady state variance
can be influenced only by the type of intervention made, that is, on
which parts of the system the interventions are made. Furthermore, by
considering explicit formulas for the steady state variance in the
original system, it may be seen that for example positive covariances
may turn negative and vice versa when making interventions.

\bigskip
\textbf{Acknowledgements.} The development of the notion of intervention for
  SDEs is joint work with my thesis advisor, Niels Richard
  Hansen, whom I also thank for valuable discussions and advice.

\bibliographystyle{amsplain}

\bibliography{full}

\end{document}